\documentclass[11pt]{amsart}
\usepackage{epsfig}
\usepackage{verbatim}
\usepackage{amssymb}
\newcommand{\eproof}{\hfill\rule{2.2mm}{3.0mm}}

\newcommand{\Proof}{\noindent {\bf Proof.~~}}

\newcommand{\HH}{{\mathcal H}}

\renewcommand{\eqref}[1]{(\ref{#1})}

\newcommand{\bfs}{{\mathbf s}}

\newtheorem{prop}{Proposition}

\newtheorem{theo}[prop]{Theorem}

\begin{document}
\baselineskip 18pt
\title{The Three Hat Problem}
\author{Brian Benson}
\email{gth858n@mail.gatech.edu}
\author{Yang Wang}
\address{School of Mathematics  \\ Georgia Institute of Technology \\
Atlanta, Georgia 30332, USA.}
\email{wang@math.gatech.edu}
\begin{abstract}
In this paper we study the Three Hat Problem which appeared in {\em Puzzle Corner} of the {\em Technology Review} magazine.  This puzzle gives a scenario in which three players wearing hats are sitting together and each hat can be seen by everyone except the player that is wearing that hat. Each player is told that all of the hats contain a positive integer and that two of the integers add to the third. In an ordered, turn-wise, modular fashion, each player truthfully states whether or not he knows his integer.  We give a strategy which allows for one of the players to solve for his integer for all possible integer configurations of the puzzle and prove it is the optimal such strategy.
\end{abstract}
\maketitle

\section{Introduction}
\setcounter{equation}{0}

Many classical puzzles involve hats. The general setting for these puzzles is a game in which
several players are each given a hat to wear.  Associated with each hat is either a color or a number. 
Each player can see the color of or number on the other players' hats, but cannot see his own.
The objective of the players is to figure out the colors or the numbers on 
their own hats. The Three Hat Problem is one of such puzzles. 

\vspace{2mm}
\noindent
{\bf The Three Hat Problem.}  Three players are each given a hat to wear. Written on each hat is
a positive integer. Any player can see the two numbers on the other players' hats, but not his own.
It is known that one of the numbers is the sum of the other two.
Individually, each player takes a turn in which he either identifies his number or passes if he cannot.  If each player fails to identify his number on his first turn, the players start the turn-wise process over again using the same order from the previous round.  The process ends when one of the players is able to correctly identify his number.
One popular example of the puzzle gives the following scenario:
\begin{itemize}
\item[] Player A:~~Pass.
\item[] Player B:~~Pass.
\item[] Player C:~~Pass.
\item[] Player A:~~My number is 50.
\end{itemize}
The question is: What are the other numbers?

There is also a more complex version of the above problem in which the players take longer to reach a solution.  It proceeds as follows:

\begin{itemize}
\item[] Player A:~~Pass.
\item[] Player B:~~Pass.
\item[] Player C:~~Pass.
\item[] Player A:~~Pass.
\item[] Player B:~~Pass.
\item[] Player C:~~Pass.
\item[] Player A:~~Pass.
\item[] Player B:~~Pass.
\item[] Player C:~~My number is 60.
\end{itemize}
Again the question is: What are the other numbers?

The most general form of the Three Hat Problem would have numbers $a, b, a+b$. In this general
setting, one may ask: (a) Will the players be able to determine their numbers, and (b) if so, how will the players 
proceed in doing so?

As far as we know, both puzzles were proposed by Donald Aucamp in
the {\em MIT Technology Review}, see 
\cite{TR03},\cite{TR04},\cite{TR06}.  Although by no means
trivial, the first puzzle is readily within grasph of most enthusiasts who have some
familiarity with these type of puzzles. The solution is Player B has 20 and Player C has 30.
To see why these two numbers work, we proceed with the following observation.  Player A on his first turn obviously doesn't know whether his
number is 50 or 10. Similarly neither Player B nor Player C can immediately figure out their numbers.
However, on his second turn Player A can reason: {\em If mine is a 10, then Player C would know
his number is either 10 or 30. If it is 10 Player B would immediately know his number is 20.
But he didn't know. So Player C should know his number is 30. Now since Player C didn't know,
my number must be 50.} With this kind of reasoning we can also rule out all other
combinations. So $[50, 20, 30]$ is the only solution to the first puzzle.
In a private communication Aucamp mentioned that he received no solution to the second
puzzle from the readers \cite{Auc-pri}. As it turns out, our study shows that the second
puzzle has eight solutions! They are 
$[25,35,60]$, $[35,25,60]$, $[42, 18, 60]$, $[18, 42, 60]$, $[10, 50, 60]$, $[50, 10, 60]$,
$[44, 16, 60]$, $[16, 44, 60]$.


The Three Hat Problem is among the more challenging hat puzzles. However, as
we shall see, many of these hat puzzles can be solved using the 
same principles and techniques as the Three Hat Problem. To illustrate, we list two classical hat puzzles.

\vspace{2mm}
\noindent
{\bf The Two Hat Problem.} Two players are each given to wear a hat with a positive integer 
written on it. Assume that the two numbers are consecutive integers.
Each player can see the other player's number but not his own. 
They take turns to either identify their numbers or pass if they cannot. 
Will the players be able to identify their numbers, and if so, what will they proceed in doing so?

\vspace{2mm}
\noindent
{\bf The Color-Hat Problem.} Several players are each given either a red or a blue
hat to wear. Each player can see all other hats but not his own. They are also told that there is
at least one red hat. The game proceeds by rounds. In each round, every player will either
identify the color of his hat or pass, but all players do so {\em simultaneously}. The game
ends when one or more players have correctly identified their colors while
none of the players responds with the incorrect color. What will happen? This puzzle takes on many popular forms, one of
which is the {\em Muddy Face Problem} analyzed in Tanaka and Tsujishita \cite{TaTs95}.

A very challenging variation of the Color-Hat Problem was due to Todd Ebert \cite{Ebe98}
and was reported
in an article in the New York Times \cite{NYT01}. In this variation, the players are
allowed to collaborate as a team and decide on a strategy before the game starts.
However, the players have only one chance to identify their colors. They win
if at least one player correctly names the color of his hat while no one is wrong. The
question is: How well can they do? What is their optimal strategy? This problem has an
interesting connection to coding theory.

Interestingly enough, each of the hat puzzles mentioned here can have a similar {\em collusion version}
that is phrased as a game of strategy. Suppose that we say the players win if at least one player makes a correct identification while none of the other players makes an incorrect identification attempt. Then each aforementioned hat puzzle can be viewed as a problem of finding the strategy for the players to win with the least
number of go-arounds. 

Although this paper is concerned with the Three Hat Problem, a main additional objective is
to show that these type of puzzles can be analyzed easily if we first treat them as games
of strategies. Once optimal strategies are found we can often easily show that the
non-collusion version and the collusion version for those games are equivalent, and therefore 
they will end in exactly the same fashion. One of the main advantages of 
presenting these puzzles as games of strategy is that we can avoid invoking the
so-called {\em super-rationality assumption} (see Hofstadter \cite{Hofs-book}), namely
each player has unlimited mental capacity to process all informations available to them, including
long chains of
reasonings such as ``I know player B knows player C knows I know player C knows ....''
Such an assumption can be confusing even to mathematicians without venturing deeply into
the realm of set theory and mathematical logic. 
The Three Hat Problem is an excellent example to illustrate this point.

\section{Optimal Strategy for the Three Hat Problem}

We now discuss a strategy for the collusion version of the Three Hat Problem. We say a strategy is 
{\em viable} if it always leads to a win for the players; in otherwords, the players need not resort to guessing at any stage.  A viable strategy is {\em optimal} if it requries the least number of turns (go-arounds) to
end the game successfully regardless what the numbers are on the three hats. Of course, while all optimal strategies must be viable, not all viable strategies are optimal. In theory it is also possible that an optimal
strategy does not exist, in which case a strategy may be the best for some configurations but no strategy is the best for all configurations.  For the Three Hat Problem, there does exist an optimal strategy which is given herein.
The optimality of the strategy is proven in the next section.

The optimal strategy we describe here is a reduction scheme involving a chain of vectors with postive integer
entries. Throughout this paper we assume that the game begins with Player A, Player B, and Player C taking turns respectively in that order.  Further, this order remains in all subsequent rounds until the game ends. The numbers $a,b,c$ for Players A, B and C respectively are represented by the vector $[a,b,c]$.
Such a vector is called a {\em Three Hat configuration}, or simply just a {\em configuration}.

Let $\HH$ denote the set of all triples $\bfs =[a,b,c]$ where $a,b,c$ are positive
integers such that the largest of which is the sum of the other two. $\HH$ represents the
set of all possible configurations of the Three Hat Problem. Define a map
$\sigma: \HH \longrightarrow \HH$ as follows: For $\bfs =[a,b,c]\in\HH$, if two of the entries are
identical then $\sigma(\bfs)=\bfs$; otherwise the largest entry is replaced by the
difference of the other two entries. For example, $\sigma([3,10, 7])=[3,4,7]$,
$\sigma([10,1, 9])=[8,1,9]$, and $\sigma([3,3, 6])=[3,3,6]$. We shall call $\bfs\in\HH$
a {\em base configuration} if $\bfs$ contains two identical entries, or equivalently
$\sigma(\bfs)=\bfs$. Note that in the base configuration, the player with the largest
number can immediately declare that his number is the sum of the other two numbers.  Although, we must note that in this situation, the aforementioned player may choose not to initially identify his number in order to obey his strategy.

Our strategy for the Three Hat Problem involves a chain of configurations for each player.
For any $\bfs\in \HH$ we otain a sequence of configurations $\bfs, \sigma(\bfs), ...
\sigma^n(\bfs)$ where $n\geq 0$ is the smallest power such that $\sigma^n(\bfs)$ is
a base configuration. For example, for $\bfs=[3,10,7]$ the sequence is
$$
   [3,10,7],[3,4,7],[3,4,1], [3,2,1], [1,2,1].
$$
We call the sequence in reverse order the {\em configuration chain} associated with $\bfs$.
So in the above example $\bfs=[3,10,7]$ the associated configuration chain is
$$
   [1,2,1], [3,2,1],[3,4,1], [3,4,7], [3, 10,7].
$$
Given a configuration,  we say that a player has
the {\em cue} if his number is the sum of the other two.
For example, for the configuration $[3,10,7]$ Player B has
the cue.

\vspace{2mm}
\noindent
{\bf Chain Reduction Strategy for the Three Hat Problem.}~~For the Three Hat Problem with 
configuration $\bfs=[a,b,c]$, let $\bfs_A=[b+c,b,c]$, $\bfs_B=[a,a+c,c]$ and $\bfs_C=[a,b,a+b]$.
These are the {\em working configurations} for Players A, B, and C respectively. 
Each player now writes down the configuration chain associated with his working configuration.
It is important to note that the chains differ only at the end. The players with the two smaller
numbers have longer chains by one configuration, which may differ for these two players. The remainder of
the chains are identical.

When the game begins, the players are assigned the first configuration in their respective
configuration chain and proceed with the following reduction scheme:

At each turn, a player looks at the configurations remaining on his configuration chain. If it 
contains only one configuration, he declares that his number is the sum of the other two numbers ending the game; otherwise, he will pass. Each player will now examine his assigned configuration (which is identical for all
the players before the game ends). If and when the player with the cue for the current configuration passes, all the players will cross out or exclude the current configuration from his chain and assign himself the next configuration
in the chain.  Otherwise, he keeps his currently assigned configuration and his chain intact.  The game continues until a player declares his number.
\eproof

\vspace{2mm}
The following two examples will facilitate the understanding of the strategy.

\vspace{2mm}
\noindent
{\bf Example 1.} The numbers for Players A, B, C are $60, 36, 24$, respectively. In this case
the working configurations are $\bfs_A=[60, 36,24]$, $\bfs_B=[60, 84, 24]$ and
$\bfs_C=[60, 36, 96]$.
The configuration chains are
$$
\begin{array}{ll}
    \text{Player A}: &  [12,12,24],[12,36,24],[60,36,24] \\
    \text{Player B}:&  [12,12,24],[12,36,24],[60,36,24],[60,84, 24] \\
    \text{Player C}: & [12,12,24],[12,36,24],[60,36,24],[60,36, 96] 
\end{array}
$$
At the start of the game, all players are assigned the configuration $[12,12,24]$.
Player A will pass, as will Player B and Player C. But Player C has the cue. So
after Player C has passed the configuration $[12,12,24]$ is crossed out by all players from
their chain. The new configuration chains are
$$
\begin{array}{ll}
    \text{Player A}: &  [12,36,24],[60,36,24] \\
    \text{Player B}: &  [12,36,24],[60,36,24],[60,84, 24] \\
    \text{Player C}: &  [12,36,24],[60,36,24],[60,36, 96] 
\end{array}
$$
All three players are now assigned the configuration $[12,36,24]$.
Player A and Player B will pass again. But since Player B has the cue, after his pass
all three players will cross out $[12,36,24]$ from their chain and assign themselves the
next configuration, which is $[60,36,24]$ for everyone. The new configuration chains are
$$
\begin{array}{ll}
    \text{Player A}: & [60,36,24] \\
    \text{Player B}: & [60,36,24],[60,84, 24] \\
    \text{Player C}: & [60,36,24],[60,36, 96] 
\end{array}
$$
It is Player C's turn and he will pass. 
Now Player A has only one configuration left on his chain, namely $[60, 36, 24]$.
So he declares his number to be the sum of the other
two numbers, which is 60. The game ends with a win for the players. 
\eproof

\vspace{2mm}
\noindent
{\bf Example 2.} The numbers for Players A, B, C are $3, 10, 7$, respectively. In this case
the working configurations are $\bfs_A=[17, 10, 7]$, $\bfs_B=[3, 10 ,7]$ and
$\bfs_C=[3, 10, 13]$. The following shows the configuration chains and the action at
each turn. Players with the cue are denoted by a $*$.
$$
\begin{array}{lcl}
\text{Player A:} & \text{Pass} & [1,2,1], [3,2,1],[3,4,1], [3,4,7], [3, 10,7], [17,10,7]\\
\text{Player B*:} & \text{Pass} & [1,2,1], [3,2,1],[3,4,1], [3,4,7], [3, 10,7] \\
\text{Player C:} & \text{Pass} & [3,2,1],[3,4,1], [3,4,7], [3, 10,7], [3, 10, 13] \\
\text{Player A*:} & \text{Pass} & [3,2,1],[3,4,1], [3,4,7], [3, 10,7], [17, 10, 7] \\
\text{Player B*:} & \text{Pass} & [3,4,1], [3,4,7], [3, 10,7] \\
\text{Player C*:} & \text{Pass} & [3,4,7], [3, 10,7], [3,10,13] \\
\text{Player A:} & \text{Pass} &  [3, 10,7], [17,10,7] \\
\text{Player B*:} & \text{\em I have 10} &  [3, 10,7].
\end{array}
$$
The game ends successfully for the players.
\eproof

Using this strategy, the player with the sum of the other two
numbers will always be the one to declare his number correctly
to end the game. This is quite easily shown. Since his chain is a
subchain of the other two players, and by the time his chain is down to
only one configuration the other players still have two. Moreover, since he
holds the cue at that stage the other players cannot reduce the chain further without
waiting for him to act. But when he does act he will declare his number. So he is always
the first to identify his number.

\section{Optimality of the Chain Reduction Strategy}

We will now prove that the above strategy is optimal for the Three Hat Problem in the sense that
no other viable strategy will be able to end the game with fewer turns for all configurations.
Before proceeding further we first notice that because $\gcd (a,b) = \gcd (a,c) = \gcd (b,c)$
the players can always divide out the numbers by the greatest common divisor of the two numbers
they see. So we may without loss of generality assume that all numbers in the three hat game are pairwise
coprime. In the coprime case the only base configurations are $[1,1,2]$, $[1,2,1]$ and $[2,1,1]$.

\begin{prop}  \label{prop-1}
     No matter what viable strategy the players use for the Three Hat Problem, the
player whose number is the sum of the other two is always the first player to declare his number.
\end{prop}
\Proof Assume that in the Three Hat Game a player declared his number on the very first turn of
the game. It is easy to see that this can happen only if we have a base
configuration and this player has the sum of the other two numbers. No other cases allow the game
to end on the very first turn without guessing.
For instance, even in the base configuration $[1,2,1]$ Player A cannot declare his number
on his first turn without guessing, for he can have both $1$ or $3$.

If the proposition is false then we have a game with configuration $[a,b,c]$ that ends 
on the $n$-th turn, with $n>1$, with a player who does not have the sum of the two numbers.
Without loss of generality, we assume that Player C declares his number to end the game, and
he does not have the sum. This gives that $c=|a-b|$; but if so, Player C must have concluded on the
$n$-th turn that his number is not $c=a+b$. This is equivalent to saying that had his number
been $c=a+b$ the game would have ended earlier, with another player declaring his number.
Therefore the strategy the players use allows them to end the 
three hat configuration $[a,b,a+b]$ in
$k<n$ turns by a player other than Player C. This player does not have
the sum of the other two numbers.

We can repeat this reasoning. In the end, we deduce that using their strategy the players 
can end a non-base configuration game in one turn by a player whose number is not the sum of the other two numbers. This is a contradiction.
\eproof

\begin{theo} \label{theo-2}
The Chain Reduction Strategy is the optimal strategy for the Three Hat Problem.
\end{theo}
\Proof For the Three Hat Problem with the configuration $[a,b,c]$ let $r([a,b,c])$
denote the number of turns needed to end the game using the Chain Reduction Strategy. We
prove that one cannot end the game in fewer turns using any other strategy.

Assume that the players are using another viable strategy such that the game ends in
$f([a,b,c])$ turns. Our objective is to show $f([a,b,c]) \geq r([a,b,c])$.
Without loss of generality, we assume that $a,b,c$ are pairwise coprime. We will prove the optimality
of the Chain Reduction Strategy by induction on $\max (a, b, c)$. 

For $\max(a,b,c)=2$ we have the base case. It is clear that
the Chain Reduction Strategy is optimal, $f([a,b,c]) \geq r([a,b,c])$.
Now assume that $f([a,b,c]) \geq r([a,b,c])$ whenever $\max(a,b,c) <M$. We now prove that
$f([a,b,c]) \geq r([a,b,c])$ if $\max (a,b,c)=M$.

We shall examine the case $a=b+c$ and $b>c$, so $a=M$. The other cases are proved in virtually 
identical fashion so we shall omit them. Note that by Proposition \ref{prop-1} the game
will end with Player A declaring his number regardless of the strategy. With this in mind
we need only to examine what happens before Player A declares his number. Clearly from
his perspective Player A knows he has either $a=b+c$ or $a=b-c$. He is not 
able to declare his number until he rules out $a=b-c$, regardless of the strategy the players
are using. Now since all strategies end with
the player with the sum declaring his number, Player A knows that if his number is
$a=b-c$ Player B will declare his number first on the $n$-th turn, where $n=
f([b-c,b,c])$. But by the $n$-th turn Player B will pass because he does not have the
sum, and after it the earliest Player A can declare his number is after Player C's pass.
Thus,
$$
	f([a,b,c]) \geq 2+ f([b-c,b,c]).
$$
Note that here we do not get equality in general because we do not assume the strategy is
optimal. By the induction hypothesis, since $\max(b-c, b, c) = b<a =M$ we have
$f([b-c,b,c]) \geq r([b-c,b,c])$, and hence $f([a,b,c]) \geq 2+ r([b-c,b,c])$.
We argue that $r([a,b,c]) = 2+r([b-c,b,c])$. This can be seen easily if we compare the 
configuration chains for $[b-c,b,c]$ and those for $[a,b,c]$. For all three players the
former is a sub-chain of the latter with one less configuration. On the $r([b-c,b,c])$-th 
turn, Player B has the cue and will pass. So $[b-c,b,c]$ is crossed out from
eveyone's chain, leaving Player A with only one configuration on his chain, namely
$[a,b,c]$. After Player C passes Player A is able to declare his number as
$a=b+c$ using the Chain Reduction Strategy.
Thus $f([a,b,c]) \geq 2+r([b-c,b,c])=r([a,b,c])$. This proves the optimality
of the Chain Reduction Strategy.
\eproof

One may wonder whether there are indeed non-optimal viable strategies for the Three Hat
Problem. One such strategy is the following: Players will note the larger of the two
numbers they see, call these $n_A$, $n_B$, and $n_C$ respectively. Unless another
player has already declared his number, Player A will pass until his $n_A$-th turn,
when he will declare his number to the the sum of the two other numbers; further, Players B and C
do likewise. This is clearly a viable strategy but by no means an optimal one.

\section{Equivalence of Collusion and No-Collusion Versions}

We now argue that under the super-rationality assumption the no-collusion
version of the Three Hat Problem will end exactly the same way as if the players are colluding using the Chain Reduction Strategy.  Specifically, we assert that if there exists an optimal strategy
then a super-rational player is able to obtain this result. Clearly, from this perspective, 
if an optimal strategy exists then the players need not collude. The super-rationality
assumption suffices to gurantee that all players will be able to find it and use it with
the knowledge that other players will do likewise. Collusion is helpful only when
there exists no single optimal strategy. This is the case when for any one strategy there
is another strategy that is better for some configurations.  If so
the players need to collude to decide on one strategy.  
Note that two strategies for the Three Hat Problem are considered to be the same if they
lead to exactly the same solution for all configurations. 
In this sense the Chain Reduction Strategy is clearly the unique optimal strategy. By the above argument we have

\begin{theo}  \label{theo-3}
     The no-collusion Three Hat Problem is equivalent to the collusion Three Hat Problem
     using the Chain Reduction Strategy.
\end{theo}

By establishing the equivalence of collusion and no-collusion versions we can also
solve the other two hat problems easily. For the Two Hat Problem, the no-collusion version
is equivalent to players using the following strategy: Each player will pass until on his
$n$-th turn, when he will declare his number to be $n+1$, 
where $n$ is the number written on the other player's hat. The game ends when one player
declares his number. For the Color-Hat Problem, the no-collusion version
is equivalent to this strategy: Players will each note how many red hats he sees. Say a
player sees $n$ red hats. He will then pass in the first $n$ rounds, but declares his hat
to be red on the $(n+1)$-th round. The games ends when some players declare their numbers.
These strategies are easily shown to be optimal by similar arguments for the Three Hat
Problem.

\noindent
{\bf Acknowledgement.} The authors would like to thank Ian Fredenburg and Don Aucamp for
very helpful discussions. This work was part of the NSF sponsored REU project 
at Georgia Institute of Technolgy in the summer of 2006 for the first author.
During the completion of this project the second author was serving as
a program director at the Division of Mathematical Sciences
of the National Science Foundation. The views expressed in this article are those of the authors,
and they do not necessarily represent the views of NSF.

\end{document}